\documentclass{amsart}

\usepackage[utf8]{inputenc}
\usepackage[T1]{fontenc}
\usepackage{float}
\usepackage{mathtools}
\usepackage{color}
\usepackage[english]{babel}
\usepackage{amsmath} 
\usepackage{amsfonts}
\usepackage{thmtools}
\usepackage{amssymb}
\usepackage{amsthm}
\usepackage{mathrsfs}
\usepackage{setspace}
\usepackage{graphicx}
\usepackage{bm}
\usepackage[all]{xy}
\usepackage{lscape}
\usepackage{latexsym}
\usepackage{epigraph}
\usepackage[Sonny]{fncychap}
\usepackage{tikz} 
\usepackage{tikz-cd}
\usetikzlibrary{matrix}
\usepackage{epsfig}
\usepackage{fancyhdr}
\usepackage{hyperref}
\usepackage{paralist}
\usepackage[mathcal]{eucal}
\usepackage{xcolor}
\usepackage{caption}
\usepackage{subcaption}
\usepackage{float}
\numberwithin{equation}{section}

\theoremstyle{cupthm}
\newtheorem{theorem}{Theorem}[section]

\newtheorem*{theorem*}{Theorem}
\theoremstyle{definition}

\newtheorem{example}[theorem]{Example}

\theoremstyle{remark}

\theoremstyle{question}
\newtheorem{question}[theorem]{Question}

\theoremstyle{definition}
\newtheorem*{definition*}{Definition}

\newtheorem*{corollary*}{Corollary}

\definecolor{darkred}{rgb}{0,0,0} 
\definecolor{darkgreen}{rgb}{0,0,0}
\definecolor{darkblue}{rgb}{0,0,0}

\begin{document}

\title{Nielsen equivalence in triangle groups}

\author{Ederson R. F. Dutra} \thanks{Author supported by FAPESP, São Paulo Research Foundation, grant 2018/08187-6.}
\address{Universidade Federal de São Carlos,  São Carlos, Brazil}
\email{edersondutra@dm.ufscar.br}

\begin{abstract} We extend the result  of \cite{Dutra} to generating pairs  of triangle groups, that is, we  show that  any   generating pair    of a triangle group  is represented by a  special  almost orbifold covering.
\end{abstract}


\maketitle
\section{Introduction}

Let  $G$ be a group and  $n\geqslant 1$. Recall that an \emph{elementary Nielsen transformation} on a $n$-tuple  $T=(g_1,\ldots, g_n)$ of elements of  $G$   is one of the following three types of moves:
\begin{enumerate}
\item[(T1)] For some $i\in \{1,\ldots, n\}$ replace $g_i$ by $g_i^{-1}$ in $T$.
\item[(T2)] For some $i\neq j$, $i,j\in \{1, \ldots, n\}$ interchange $g_i$ and $g_j$ in $T$. 
\item[(T3)] For some $i\neq j$, $i,j\in \{1, \ldots, n\}$ replace $g_i$ by $g_ig_j$ in $T$.
\end{enumerate}
We  say that two $n$-tuples $T$ and $T'$   are \emph{Nielsen equivalent}  if there exists a finite sequence of $n$-tuples 
$$T=T_0, T_1, \ldots, T_{k-1}, T_k=T'$$ 
such that $T_i$ is obtained from $T_{i-1}$ by an elementary Nielsen transformation for $1\leqslant i\leqslant k$. Nielsen equivalence clearly defines an equivalence relation for tuples of elements of $G$. 

\smallskip
 
Nielsen equivalence in Fuchsian groups  have been extensively studied by many authors, see  \cite{Fine, Lustig1, Lustig2, Lustig3, Lustig4, P,Rosenberger1,  R, Zieschang} for example. However the techniques deployed so far are mainly algebraic such as normal form and K-theoretic arguments.	In \cite{Dutra} we give a new geometric insight to the problem by showing that  any generating tuple of the fundamental group of a sufficiently large $2$-orbifold is represented by a special almost orbifold covering.  We still do not know if the same holds  for non sufficiently large orbifolds. The interesting case are orbifolds  with three cone points whose underlying surface is a sphere. The purpose of  this note  is to show that   generating pairs of the fundamental group of such orbifolds also have this geometric description.
\begin{theorem}{\label{mainthm}}
Let  $\mathcal{O}=S^2(p_1, p_2, p_3)$   and  
$$G=\pi_1^o(\mathcal{O})\cong \Delta(p_1, p_2, p_3)=\langle s_1, s_2, s_3 \ | \ s_1^{p_1}=s_2^{p_2}=s_3^{p_3}=s_1s_2s_3=1\rangle .$$ 

Then for any generating pair   $T$  of  $G$  there is a special almost orbifold covering  $\eta':\mathcal{O}'\rightarrow \mathcal{O}$ and a generating pair  $T' $ of $\pi_1^o(\mathcal{O}')$ such that $\eta_{\ast}'(T') $ and $T$ are Nielsen   equivalent.   
\end{theorem}

A natural question that arises  is  whether  this holds  for arbitrary generating tuples, that is: 
\begin{question}
Are all  generating tuples of a triangle group also  represented by a (special)  almost orbifold coverings?
\end{question}

The proof of Theorem~\ref{mainthm} relies heavily on the  description of all generating pairs of two-generated Fuchsian groups given  by Fine-Knapp-Matelsky-Purzitsky-Rosenberger \cite{Fine, Knapp,Matelski, P,R}.  
  

\section{Orbifolds and almost orbifold coverings}{\label{secorbifold}}

In this section we give a quick review  about  orbifolds, orbifold fundamental groups and orbifold coverings.  More details about orbifolds can be found in the beautiful  article of P. Scott~\cite{Scott}.

A  $2$-\emph{orbifold}  $\mathcal{O}$ is a pair $(F, p)$ where $F$ is a  connected  surface, called the \emph{underlying surface} of $\mathcal{O}$, and $p:F\rightarrow \mathbb{N}$ is a function such that 
$$\Sigma(\mathcal{O}):=\{x\in F \ | \ p(x)\geqslant 2\}$$ 
is finite and contained in the interior   of $F$. The number $p(x)$ is called the \emph{order} of   $x$ and the points of $\Sigma(\mathcal{O})$   are called \emph{cone points}.  An  orbifold is said to be \emph{compact} (resp. \emph{closed}) if its underlying surface is  compact (resp. closed).  

A compact orbifold   $\mathcal{O}=(F, p)$  with $\Sigma(\mathcal{O})=\{x_1, \ldots, x_r\}$  will also be denoted by   $F(p_1, \ldots, p_r)$   where   $p(x_i)=p_i$ for $i=1,\ldots, r$.

\smallskip 

The \emph{fundamental group} $\pi_1^{o}(\mathcal{O})$ of  $\mathcal{O}=(F, p)$  is defined in terms of the fundamental group of the underlying surface $F$  as follows.  Let $x_1, x_2, \ldots$ be the cone points of $\mathcal{O}$ and for each $i=1, 2 , \ldots$ let  $D_i\subseteq \text{Int} (F)$ be a disk centered at $x_i$ such that $ {D}_i\cap {D}_j= \emptyset$  for $i\neq j$.   Set   
$$ S_{\mathcal{O}} = {F -  \underset{i=1,2,\ldots}\cup \text{Int}(D_i )}.$$  
Then $\pi_1^{o}(\mathcal{O})$ is the  group obtained from $\pi_1 (S_{\mathcal{O}})$ by adding the relations $s_i^{p_i}=1$ for all $i=1,2,\ldots$,   where $s_i$ is the homotopy class represented by $\partial D_i$ and $p_i$ is the order of  $x_i$. For example the fundamental group of   $\mathcal{O}=S^2(p_1, p_2, p_3)$  has presentation   
$$\pi_1^o(\mathcal{O})=\langle s_1, s_2, s_3 \ | \  s_1^{p_1}=s_2^{p_2}=s_3^{p_3}=s_1s_2s_3=1 \rangle.$$
\begin{figure}[h!]
\begin{center}
\includegraphics[scale=0.8]{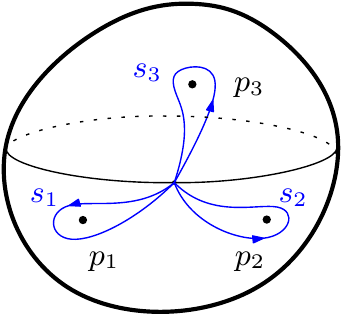}
\end{center}
\caption{The generator $s_i$ is represented by a simple closed curve that goes once around the cone point $x_i$.}\label{fig:almost1}
\end{figure}

In order to define almost orbifold coverings we first  recall the definition of an orbifold covering. Let $\mathcal{O}=(F,p)$ and $\mathcal{O}'=(F', p')$ be orbifolds. An \emph{orbifold covering} $\eta:\mathcal{O}'\rightarrow \mathcal{O}$ is a continuous surjective map $\eta:F'\rightarrow F$ 
with the following properties:
\begin{enumerate}
\item For each point $y\in F'$, the order of $y$ divides the order of $\eta(y)$.

\item For each point $x\in \text{Int}(F)$ the set $\eta^{-1}(x)\subseteq F'$ is discrete and, over a small disk in $F$ centered at $x$,  the map $\eta$ is equivalent to the map 
$$(z, y)\in \mathbb{D}^2 \times \eta^{-1}(x) \mapsto e^{( \frac{2\pi p'(y)}{p(x)}i ) } z\in \mathbb{D}^2$$ 
and $x$ corresponds to $0$ in $\mathbb{D}^2$.

\end{enumerate} 
Note that the map $\eta|_{F'- \eta^{-1}(\Sigma(\mathcal{O}))}:F'- \eta^{-1}(\Sigma(\mathcal{O}))\rightarrow F-\Sigma(\mathcal{O})$ is a genuine covering.  The \emph{degree}  of  $\eta:\mathcal{O}'\rightarrow\mathcal{O}$ is defined as  the degree of $\eta|_{F'- \eta^{-1}(\Sigma(\mathcal{O}))}$. It is not hard to see that an orbifold covering $\eta:\mathcal{O}'\rightarrow \mathcal{O}$ induces a monomorphism $\eta_{\ast}:\pi_1^o(\mathcal{O}')\rightarrow \pi_1^o(\mathcal{O}).$ 
Conversely, for any subgroup $H\leq \pi_1^o(\mathcal{O})$ there is an orbifold covering $\eta:\mathcal{O}_{H}\rightarrow  \mathcal{O}$ such that $\eta_{\ast}(\pi_1^o(\mathcal{O}_{H}))=H$.

\begin{definition*}[Almost orbifold covering]{\label{def:almost}}
Let $\mathcal{O}'=(F', p')$ and $\mathcal{O}=(F, p)$ be compact  orbifolds.  An \emph{almost orbifold covering} $\eta:\mathcal{O}'\rightarrow \mathcal{O}$ is a (non necessarily surjective) continuous map $\eta:F'\rightarrow F$ having the following properties:
\begin{enumerate}
\item[(C1)] For each $y\in F'$ the order of $y$ divides the order  of $\eta(y)$. 

\item[(C2)] There is a point $x\in\text{Int}(F)$ of order $p\geqslant 1$,  called the \emph{exceptional point}, and a disk ${D}\subseteq\text{Int}(F)$ centered at $x\in F$, called the \emph{exceptional disk}, with 
$$(D-\{x\})\cap \Sigma(\mathcal{O})=\emptyset$$  
such that $\eta$ restricted to $F'-\eta^{-1}(\text{Int}(D))$ defines  an orbifold covering of finite degree between the compact $2$-orbifolds  
$$\mathcal{Q}':=(F'-\eta^{-1}(\text{Int}(D)), p'|_{F'-\eta^{-1}(\text{Int}(D))})\subseteq \mathcal{O}'$$ and $$\mathcal{Q}:=(F-\text{Int}(D), p|_{F-\text{Int}(D)})\subseteq\mathcal{O}.$$

\item[(C3)] $\eta^{-1}(D)=D_1\sqcup D_2\sqcup\ldots\sqcup D_t\sqcup C_1\sqcup \ldots \sqcup C_u$, where $t\geqslant 0$, $u\geqslant 1$  and  
\begin{enumerate}
\item[(C3.a)]  $C_1, \ldots, C_u\subseteq \partial F'$ are boundary components of $\mathcal{O}'$, called the \emph{exceptional boundary components} of $\mathcal{O}'$.

\item[(C3.b)] Each $D_j\subseteq \text{Int}(F')$ ($1\leqslant j\leqslant t$) is a disk and $\eta|_{D_j}:D_j\rightarrow D$ is equivalent to the map
$$z\in \mathbb{D}^2\longmapsto e^{\frac{2\pi p}{q}i}  z\in \mathbb{D}^2$$
and  $x\in D$ corresponds to $0$ in $\mathbb{D}^2$, where  $q$ is the order of the  point $\eta^{-1}(x)\cap D_j$.  
\end{enumerate}

\end{enumerate}
The degree of $\eta:\mathcal{O}'\rightarrow \mathcal{O}$ is defined as the degree of the orbifold covering  $$\eta|_{\mathcal{Q}'}:\mathcal{Q}'\rightarrow \mathcal{Q}.$$  
We further call an almost orbifold covering   $\eta :\mathcal{O}'\rightarrow \mathcal{O}$   \emph{special} if $u=1$ (i.e.   $\mathcal{Q}'$  has a single exceptional boundary component) and  the degree of the map  $\eta|_{C_1}:C_1\rightarrow \partial D$ is at most  the order $p$ of the exceptional point $x$.
\end{definition*}

\begin{example}{\label{ex:almost}}
Let  $\mathcal{O}=(F,p)$ be a closed orbifold and  $x$ an arbitrary  point of $F$. Let further   $D\subseteq F$ be a small disk centered at $x$  such that $D$ contains no cone points except possibly  $x$. Denote by $\mathcal{O}_0$ the orbifold $( F-\text{Int}(D), p|_{F-\text{Int}(D)})$ and by $j$ the inclusion map $ F-\text{Int}(D) \hookrightarrow F$. 

Any orbifold covering  $\eta':\mathcal{O}'\rightarrow \mathcal{O}_0$    of finite degree induces an almost orbifold covering  of $\mathcal{O}$,  namely the map  $\eta:= j\circ \eta':F'\rightarrow F.$
 Note that the  exceptional point is   $x$,  the exceptional disk is $D$ and 
$$\eta^{-1}(D)=\partial \mathcal{O}'= C_1\sqcup \ldots \sqcup C_u$$  are  the exceptional boundary components of $\eta$.  When   the degree of $\eta'$ is one, an hence the degree of $\eta$ is also   one,  then  $\mathcal{O}'=\mathcal{O}_0$.  In this case   we say that $\eta$   is   trivial.	
\end{example}


\section{Proof of Theorem~\ref{mainthm}}
Throughout  this section $\mathcal{O}$ will denote the orbifold $S^2(p_1,p_2,p_3)$   and $G$ its fundamental  group. In \cite{Matelski} (see also Knapp \cite{Knapp} for the particular case of  generating pairs consisting  of elliptic elements and \cite{Fine} for a more general discussion)  Matelski  implicitly shows that   if  $T$ is  a generating pair of $G$ (the arguments in \cite{Matelski} work also  for non-hyperbolic triangle groups)  then  $T$ is equivalent to $T'$ where one of the following holds:
\begin{enumerate}
\item $T'=(s_{i_1}^{\nu_1},s_{i_2}^{\nu_2})$ with  $1\leqslant i_1<i_2\leqslant 3$ and $(p_{i_j}, \nu_j)=1$ for $j=1,2$.

\item $p_1=2$,  $p_3\geqslant 3$ odd and  $T'=( s_2^{\nu}, s_1 s_2^{\nu'}  s_1^{-1})= ( s_2^{\nu}, s_1  s_2^{\nu'} s_1 )$  with  $\nu $ and $\nu'$ coprime with $p_2$. 
 
\item  $p_1=2$, $p_2=3$  and $p_3\geqslant 3$ is odd,   and   $T'=(s_1, s_2^{-1}   \cdot  s_3^{\nu}  \cdot  s_2)$ with $(\nu, p_3)=1$.

\item $p_1=2$, $p_2=3$ and $p_3\geqslant 4$  coprime with  $3$, and   $$T'=( s_1s_2s_1 \cdot  s_3^{\nu} \cdot (s_1s_2s_1)^{-1}, s_2) $$   with $\nu$ coprime with  $p_3$.

\item $p_1=2$, $p_2=3$ and $p_3\geqslant 5$  coprime with  $4$, and   
$$T'=( (s_2s_1)^2 \cdot  s_3^{\nu}\cdot    (s_2s_1)^{-2}, s_3^{\nu'})$$ 
with $ \nu$ and $\nu'$ coprime with  $p_3$.

\item  $(p_1,p_2,p_3)=(2,3,5)$, and one of the following holds:
\begin{enumerate}
\item[(6.a)]  $T'=( s_1, (s_2^2s_1)^2\cdot  s_2 \cdot  (s_2^2s_1)^{-2})$. 
 
\item[(6.b)]  $T'=(s_2\cdot s_1 \cdot s_2^{-1}, (s_1s_2)^2 s_1\cdot   s_3^{\nu}\cdot   s_1^{-1} (s_1s_2)^{-1})$ with  $\nu \in \{1, 2\}$.

 \item[(6.c)]  $T'=( s_2s_1 \cdot s_2 \cdot  (s_2s_1)^{-1}, s_3^3  s_1\cdot    s_3^{\nu} \cdot   (s_3^3s_1)^{-1}) $ with $\nu\in \{1, 2\}$.

 \item[(6.d)] $T'=(s_2  s_1 \cdot  s_2  \cdot (s_2s_1)^{-1}, (s_1s_2)^2 s_1\cdot   s_3^{\nu}\cdot    s_1^{-1} (s_1s_2)^{-1})$. 
\end{enumerate}

\item $(p_1,p_2,p_3)=(2,3,7)$, and  one of the following holds:   
\begin{enumerate}
\item[(7.a)] $T'=(s_2 ,  s_3^{-3} s_1 \cdot    s_3^{\nu} \cdot  s_1    s_3^3)$ with $ \nu$ coprime with $7$

\item[(7.b)] $T'= ((s_1s_2^2s_1s_2)^2s_2s_1s_2, (s_2(s_2s_1)^2)^2 s_2s_1) $.
\end{enumerate}

\item $p_1=p_2=3$ and $p_3\geqslant 4$  coprime with $3$,  and   $$T'= (s_1s_2^2, s_2^2s_1).$$

\item $p_1=2$, $p_2=4$ and $p_3\geqslant 7$ odd, and 
$$T'=( s_1s_2^2,  s_2^3  s_1  s_2^3).$$

\item $p_1=2$, $p_2=3$ and $p_3\geqslant 7$ coprime with  $6$, and  $$T'=(s_1s_2s_1s_2^2, s_2^2s_1s_2s_1).$$
\end{enumerate}
Note that the list (1)-(10) above does not give a classification of generating pairs, that is, pairs from distinct items might be Nielsen equivalent.  For example the   generating  pair $T'=(s_2^{\nu}, s_1s_2^{\nu'}s_1^{-1})$ given in (2) with $\nu=\nu'=1$   is Nielsen equivalent to the standard pair  $(s_1, s_3^2)$.  For a complete classification of generating pairs of triangle groups see \cite{Fine}.  
 
For  each   generating pair given in  (1)-(10) we will  describe a  special almost orbifold covering that represents $T$.  First note that standard generating tuples, that is, those given in   (1),  are represented by trivial almost orbifold coverings, see Example~\ref{ex:almost} and  \cite[Example 1.9]{Dutra}.  Thus  we are left with the generating pairs given in  items  (2)-(10). 

\smallskip

\begin{figure}[h!]
\begin{center}
\includegraphics[scale=0.8]{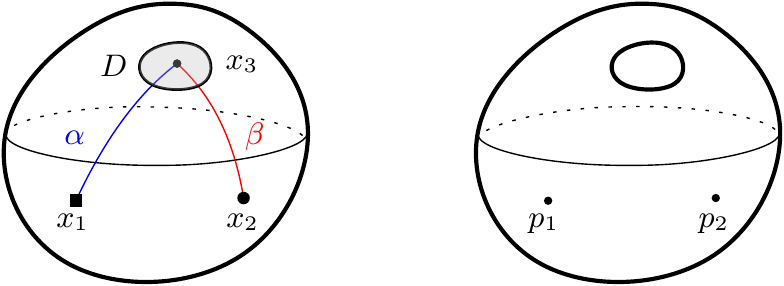}
\end{center}
\caption{ The  arcs  $\alpha $ and $ \beta$ in $S^2$  and the orbifold $\mathcal{O}_0=D_0(p_1,p_2)\subseteq \mathcal{O}$.}\label{fig:case0}
\end{figure}

Let $D\subseteq S^2$ be a disk centered at the   cone point  $x_3$  that does not contain $x_1$ and $x_2$ as shown in Fig.~\ref{fig:case0}. Denote the orbifold obtained from $\mathcal{O}$ by removing the interior of  $D$  by $\mathcal{O}_0$. Then $\mathcal{O}_0$  is a disk  $D_0:=S^2-\text{Int}(D)$ with two cone points $x_1$ and $x_2$ of order $p_1$ and $p_2$ respectively. Consequently  $\pi_1^o(\mathcal{O}_0)$ has the following presentation
$$ \langle S_1 , S_2 \ | \ S_1^{p_1}=S_2^{p_2}=1\rangle.$$  
The  inclusion map $i:D_0\rightarrow S^2$  induces a surjective  homomorphism   $$i_{\ast} :\pi_1^{o}(\mathcal{O}_0)\twoheadrightarrow \pi_1^o(\mathcal{O})$$
 that sends $S_i$ onto $s_i$ ($i=1,2$)  and $S_1S_2$ onto $s_3^{-1}$. 

Let further  $\alpha$ and $\beta$    be arcs in $S^2$  with $\partial \alpha =\{x_1, x_3\} $ and  $\partial \beta= \{x_2, x_3\}$  as shown in   Fig.~\ref{fig:case0}.   Let further $\alpha_0$ and  $\beta_0$ be the arcs $D_0\cap \alpha $ and $ D_0\cap \beta$  in $D_0$. We also  label the cone point   $x_1$ (of order $p_1$) by a small square an the cone point $x_2$ (of order $p_2$) by a small disk. The  arcs $\alpha$, $\alpha_0$, $\beta$ and $\beta_0$   and  the labels of the points  $x_1$ and $x_2$  will help us to visualize the defined orbifold maps. 

In all cases we will describe an almost orbifold covering $\mathcal{O}'\rightarrow \mathcal{O}$  where the   orbifold $\mathcal{O}'$  is either a disk with two cone points   or a torus with no cone points that has an open disk removed. Moreover, in all cases   $x_3$ will be  the exceptional point and $D$  the exceptional disk.

 \smallskip 

\noindent \textbf{Case (2).} Recall that $p_1=2$ and that  $p_3\geqslant 3$ is odd. The pair $T'$ is equal to $(s_2^{\nu} , s_1s_2^{\nu'}s_1^{-1})$ with $(\nu, p_2)=(\nu', p_2)=1$. 

\begin{figure}[h]
\begin{center}
\includegraphics[scale=0.8]{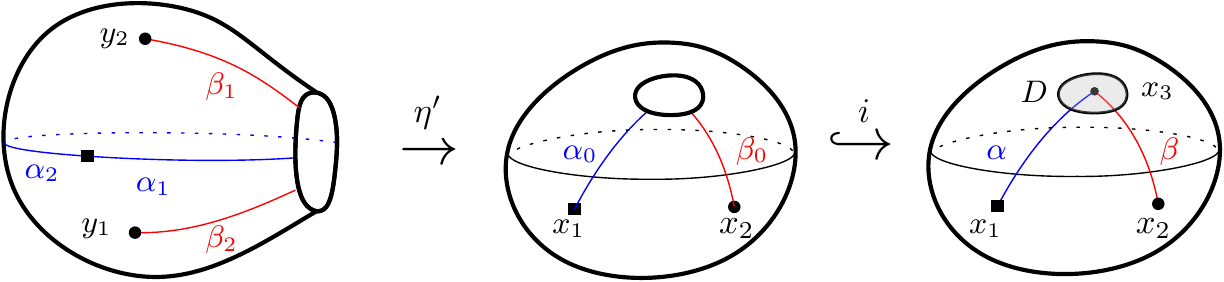}
\end{center}
\caption{ $\eta=\eta'\circ i:D^2\rightarrow S^2$ defines a special almost orbifold covering of degree $2$. }\label{fig:case2}
\end{figure}

Consider the surjecitve  map $\eta':D^2\twoheadrightarrow D_0$ described in Fig.~\ref{fig:case2}, that is, $\eta'(\alpha_1)=\eta'(\alpha_2)=\alpha_0$, $\eta'(\beta_1)=\eta'(\beta_2)=\beta_0$ and $\eta'$  maps each component of $$D^2-\alpha_1\cup \alpha_2\cup \beta_1\cup \beta_2$$ homeomorphically onto $D_0-\alpha_0\cup   \beta_0$.
 
 Denote by  $\mathcal{O}'$   the orbifold with underlying surface $ D^2$   and with cone points $y_1$ and $y_2$ both of order $p_2$. Thus $\mathcal{O}'=D^2(p_2,p_2)$.  Then  it is not hard to see that $\eta'$ defines an orbifold covering of degree $2$  from $\mathcal{O}'$ onto $\mathcal{O}_0\subseteq \mathcal{O}$  that corresponds to the subgroup $\langle S_2 , S_1S_2S_1^{-1}\rangle$ of $\pi_1^o(\mathcal{O}_0)$. As $\eta_{\ast}':\pi_1^o(\mathcal{O}')\rightarrow \pi_1^o(\mathcal{O}_0)$ is 1-1 we  may assume that $$\pi_1^o(\mathcal{O}')=\langle S_2 , S_1S_2S_1^{-1}\rangle.$$ 

It follows from Example~\ref{ex:almost} that the map $\eta:=i\circ \eta':D^2\rightarrow S^2$ defines an almost orbifold covering from $\mathcal{O}'$ to $\mathcal{O}$. As the degree of $\eta|_C:C\rightarrow \partial D$ is equal to two, which is strictly smaller than $p_3$, we conclude that $\eta$ is   special. Moreover, 
$$ \eta_{\ast}(S_2^{\nu})=i_{\ast}(S_2^{\nu})= s_2^{\nu} \ \  \text{ and } \ \   \eta_{\ast}(S_1S_2^{\nu'}S_1^{-1})=i_{\ast}(S_1S_2^{\nu'}S_1^{-1})=s_1s_2^{\nu'}s_1^{-1}$$ 
which means that $\eta$ represents the generating pair $T'$.
  
\smallskip

\begin{figure}[h!]
\begin{center}
\includegraphics[scale=0.8]{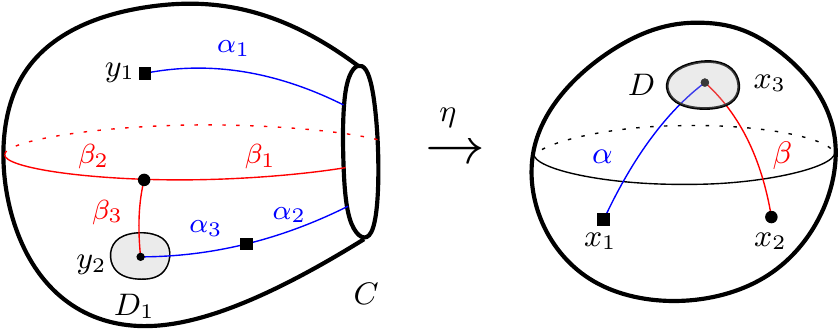}
\end{center}
\caption{ Case \textbf{(3)}.  $\eta:D^2\rightarrow S^2$ defines a special almost orbifold covering of degree $3$.  }\label{fig:case3}
\end{figure}

\noindent\textbf{Case  (3).} Recall that in this case we have $p_1=2$, $p_2=3$ and $p_3\geqslant 3$ odd, and  $T'=(s_1, s_2^{-1}    s_3^{\nu}   s_2)$ with $(\nu, p_3)=1$.  
Let $\eta:D^2\rightarrow S^2$ be the map  described in Fig.~\ref{fig:case3}, that is,   $\eta(\alpha_1)=\eta(\alpha_2)=\alpha_0$, $\eta(\alpha_3)=\alpha$, $\eta (\beta_1)=\eta (\beta_2)=\beta_0$,  $\eta (\beta_3)=\beta$ and $\eta $ maps each component of 
$$D^2-\text{Int}(D_1)\cup \alpha_1\cup\alpha_2\alpha_3\cup\beta_1\cup\beta_2\cup \beta_3$$
homeomorphically onto $D_0-\alpha_0\cup \beta_0$ and the disk  $D_1$ homeomorphically onto $D$.

Denote by $\mathcal{O}'$   the orifold  with underlying surface $D^2$ and with two cone pints $y_1$  and $y_2$ of order $p_1=2$ and $p_3$ respectively, i.e. $\mathcal{O}'=D^2(p_1, p_3)$.   Then  $\eta$  defines an almost orbifold covering of degree $3$ from $\mathcal{O}'$ to $\mathcal{O}$.

By the description of $\eta$ we see that   $\eta^{-1}(D)=D_1 \sqcup  C$.  The restriction of $\eta$ to the disk $D_1$ defines an orbifold covering of degree $1$ from the orbifold $\mathcal{Q}'=D_1(p_2)$ (a disk with a single cone point, $y_2$ ,  order $p_3$)   onto  $\mathcal{Q}=D(p_3)$ (a disk with a single cone point, $x_3$ , of order $p_3$).  Moreover,   the exceptional boundary component $C$ is mapped  by $\eta$ onto $\partial D$ with degree $2$ which implies that   $\eta$ is a special almost orbifold covering. 

Put $F:= D^2-\text{Int}(D_1)$. It is easy to see that  the map 
$$\eta':=\eta|_{F}:F:=D^2-\text{Int}(D_1)\twoheadrightarrow D_0$$  defines an orbifold covering of degree $3$ from the orbifold $\mathcal{O}_0'=F(p_1)$ (an annulus  with a single cone point $y_1$ of order $p_1=2$) onto the orbifold $\mathcal{O}_0\subseteq \mathcal{O}$ and that   $\eta'$ corresponds to the subgroup  $\langle  S_1 , S_2^{-1}\cdot  S_1S_2\cdot  S_2\rangle$ of $\pi_1^o(\mathcal{O}_0)$. Thus we may assume that   
$$\pi_1^o(\mathcal{O}_0')=\langle S_1 , S_2^{-1}\cdot S_1S_2\cdot S_2\rangle.$$
The inclusion map $j:F\hookrightarrow D^2$  induces a surjective homomorphism $$j_{\ast}:\pi_1^o(\mathcal{O}_0')\rightarrow \pi_1^o(\mathcal{O}').$$ Consequently $\pi_1^o(\mathcal{O}')$ is generated by  
$$t_1:=j_{\ast}(S_1) \ \ \text{ and } \ \ t_2:=j_{\ast}(S_2^{-1}\cdot S_1S_2 \cdot S_2).$$     
As  $\eta_{\ast}\circ j_{\ast}=i_{\ast}\circ \eta_{\ast}'$ and   $i_{\ast}(S_1S_2)=s_3$,  we see that $\eta_{\ast}$ sends the generating pair $(t_1, t_2^{-\nu})$  of $\pi_1^o(\mathcal{O}')$   onto the pair $T'$, and hence  the special almost orbifold covering $\eta:\mathcal{O}'\rightarrow \mathcal{O}$ represents the generating pair $T'$.

\smallskip 

\begin{figure}[h!]
\begin{center}
\includegraphics[scale=0.8]{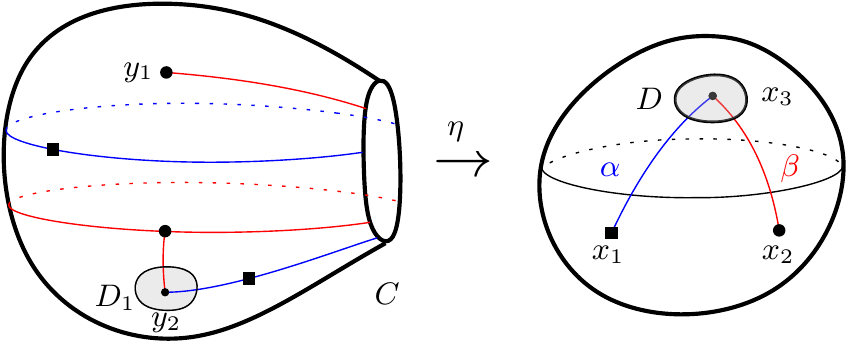}
\end{center}
\caption{ Case \textbf{(4)}.  $\eta:D^2\rightarrow S^2$ defines a special almost orbifold covering of degree $4$. }\label{fig:case4}
\end{figure}

\noindent \textbf{Case (4).}  In this case we have $p_1=2$, $p_2=3$ and $p_3\geqslant 4$ with $(p_3, 3)=1$, and 
$$T'=( s_1s_2s_1 \cdot   s_3^{\nu} \cdot   (s_1s_2s_1)^{-1}, s_2) $$ 
with $(\nu, p_3)=1$. Consider the map   $\eta:D^2\rightarrow S^2$   described in Fig.~\ref{fig:case4} (its precise description is given as in   case (3)) and let $\mathcal{O}'$ be the orifold  with underlying surface $D^2$ and with two cone points $y_1$  and $y_2$ of order $p_2=3$ and $p_3$ respectively, i.e. $\mathcal{O}'=D^2(p_2, p_3)$.   Then $\eta$  defines an almost orbifold covering of degree $4$ from $\mathcal{O}'$ to $\mathcal{O}$. 

The description of $\eta$  implies  that   $\eta^{-1}(D)=D_1  \sqcup  C$.  The restriction of $\eta$ to the disk $D_1$ defines an orbifold covering of degree $1$ from the orbifold $\mathcal{Q}'=D_1(p_2)$ (a disk with a single cone point $y_2$  order $p_3$)   onto  $\mathcal{Q}=D(p_3)$ (a disk with a single cone point $x_3$  of order $p_3$).    Moreover,   the exceptional boundary component $C$ is mapped  by $\eta$ onto $\partial D$ with degree $3$,  and so $\eta$ is a special almost orbifold covering.

The restriction $\eta'$  of $\eta$ to the  surface  $F:=D^2-\text{Int}(D_1) \subseteq D^2$ defines an orbifold covering of degree $4$ from the orbifold $\mathcal{O}_0'=F(p_2)$ (an annulus  with a single cone point $y_1$ of order $p_2=3$) onto the orbifold $\mathcal{O}_0\subseteq \mathcal{O}$ such that 
$$\eta_{\ast}'(\pi_1^o(\mathcal{O}_0'))= \langle S_1S_2S_1   \cdot S_1S_2  \cdot S_1S_2^2S_1
, S_2\rangle. $$ 
Thus we may assume that  $\pi_1^o(\mathcal{O}_0')=\langle S_1S_2S_1   \cdot S_1S_2  \cdot S_1S_2^2S_1 , S_2\rangle.$ 

Now the inclusion map $j:F\hookrightarrow D^2$  induces a surjective homomorphism $j_{\ast}$ from $\pi_1^o(\mathcal{O}_0')$ onto $\pi_1^o(\mathcal{O}').$ Consequently $\pi_1^o(\mathcal{O}')$ is generated by  $$t_1:=j_{\ast}(S_1S_2S_1   \cdot S_1S_2  \cdot S_1S_2^2S_1) \ \ \text{ and } \ \  t_2:=j_{\ast}(S_2). $$ 
Since $\eta_{\ast}\circ j_{\ast}=i_{\ast}\circ \eta_{\ast}'$ and since $i_{\ast}(S_1S_2)=s_3^{-1}$,  we see that $\eta_{\ast}$ sends the generating pair $(t_1^{-\nu}, t_2)$  of $\pi_1^o(\mathcal{O}')$   onto the pair $T'=( s_1s_2s_1   s_3^{\nu}  (s_1s_2s_1)^{-1}, s_2)$.

\smallskip

\noindent \textbf{Case (5).} Recall that $p_1=2$, $p_2=3$ and $p_3\geqslant 5$ with $(p_3, 4)=1$, and that   $T'$ is equal to $( (s_2s_1)^2\cdot   s_3^{\nu} \cdot    (s_2s_1)^{-2}, s_3^{\nu'})$ with $\nu$ and $\nu'$ coprime with $p_3$. 

\begin{figure}[h!]
\begin{center}
\includegraphics[scale=0.8]{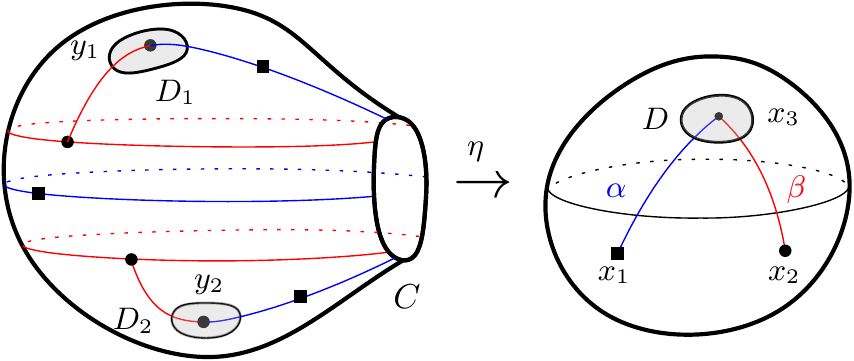}
\end{center}
\caption{ Case \textbf{(5)} $\eta:D^2\rightarrow S^2$ defines a special almost orbifold covering of degree $6$.}\label{fig:case5}
\end{figure}

Consider the  map $\eta:D^2\rightarrow S^2$  described in Fig.~\ref{fig:case5} (again, its precise definition is given as in case (3))  and let $\mathcal{O}'$ be the orifold  with underlying surface $D^2$ and with two cone pints $y_1$  and $y_2$  both of  order $p_3$, i.e. $\mathcal{O}'=D^2(p_3, p_3)$.   Then it is easy to see that  $\eta$  defines an almost orbifold covering of degree $6$ from $\mathcal{O}'$ to $\mathcal{O}$. 

Put $F:=D^2-\text{Int}(D_1)\cup \text{Int}(D_2)$.   It is not hard to see that  the map   $$\eta':=\eta|_{F}:F  \rightarrow D_0$$ defines an orbifold covering of degree $6$ from the orbifold $\mathcal{O}_0'=F$ (a pair of pants without cone points)  onto the orbifold $\mathcal{O}_0\subseteq \mathcal{O}$.  Moreover,  
$$\eta_{\ast}'(\pi_1^o(\mathcal{O}_0'))= \langle    (S_2S_1)^2 \cdot S_2^2S_1 \cdot     (S_2S_1)^{-2}, S_2^2S_1 )  \rangle.  $$ 
Thus  we can assume that    $\pi_1^o(\mathcal{O}_0')=\langle    (S_2S_1)^2 \cdot S_2^2S_1 \cdot     (S_2S_1)^{-2}, S_2^2S_1 )  \rangle$.

The inclusion map $ j:F\hookrightarrow D^2$  induces  a  surjective homomorphism $j_{\ast}$ from $\pi_1^o(\mathcal{O}_0')$ onto $ \pi_1^o(\mathcal{O}')$.  Consequently $\pi_1^o(\mathcal{O}')$ is generated by  
$$t_1:=j_{\ast}((S_2S_1)^2 \cdot S_2^2S_1 \cdot     (S_2S_1)^{-2})  \ \ \text{ and } \ \  t_2:=j_{\ast}(S_2^2S_1).$$ 
As $\eta_{\ast}\circ j_{\ast}=i_{\ast}\circ \eta_{\ast}'$ and since $i_{\ast}$ maps the element $(S_2^2S_1)$ onto $ s_3$,  we see that  the   pair $(t_1^{-\nu}, t_2)$  is mapped by $\eta_{\ast}$    onto the pair $T'$.

\smallskip

\begin{figure}[h!]
\begin{center}
\includegraphics[scale=1]{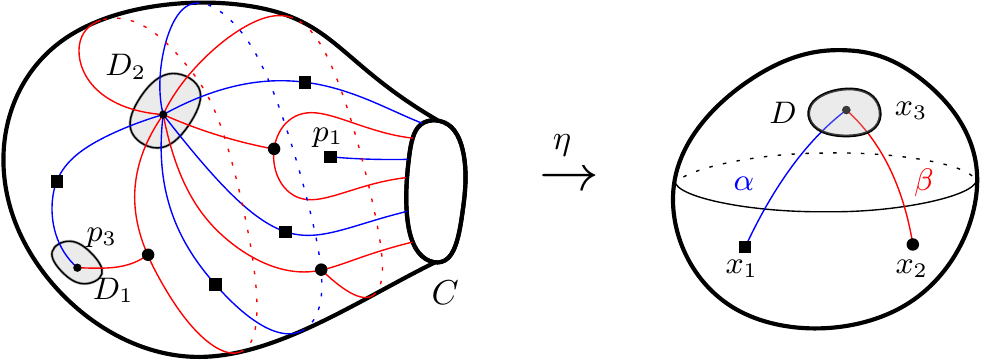}
\end{center}
\caption{ Case \textbf{(6.b)}.  $\eta$ has degree $6$ and $\mathcal{O}'$ is equal to $D^2(p_1, p_3)$.}\label{fig:case_6}
\end{figure}

\noindent \textbf{Case (6).}  $(p_1,p_2,p_3)=(2,3,5)$. It is not hard to see that    all generating pairs given in (6.a)-(6.d) are equivalent to those given in (1)-(2).  For example, suppose that    (6.b) holds, that is,     $T'=(s_2s_1s_2^{-1}, (s_1s_2)^2s_1 \cdot s_3^{\nu} \cdot s_1(s_1s_2)^{-2})$ with $\nu \in \{1,2\}$. If  $\nu =1$  then we have 
$$ T'\sim_{NE} (s_1 , s_3^2)$$ 
 which is a standard generating pair,  and if $ \nu=2$ then  we have
$$T'\sim_{NE} (s_3^2 , s_1 s_3^2  s_1^{-1})$$ 
which is case (2.b) with $p_2=3$ and $p_3=5$. We point out that for all cases (6.a)-(6.d)  it is also possible to construct a special almost orbifold covering $\eta:\mathcal{O}'\rightarrow \mathcal{O}$   and find a   generating pair $T_{\mathcal{O'}}$    of $\pi_1^o(\mathcal{O}')$ that is mapped by $\eta_{\ast}$ onto $T'$.  If we take  case (6.b)   for example,  then  the almost orbifold covering  $\eta$ described in Figure~\ref{fig:case_6} represents $\tau'$. 

\smallskip  
 
\noindent \textbf{Case (7).} $(p_1, p_2, p_3)=(2,3,7)$.  We start with the generating par  given in (7.a), that is, we define  a special  almost orbifold covering  that represents  
$$T'=(s_2 ,  s_3^{-3} s_1 \cdot    s_3^{\nu} \cdot  s_1    s_3^3)$$ 
where $(\nu, 7)=1$.

\begin{figure}[h!]
\begin{center}
\includegraphics[scale=0.8]{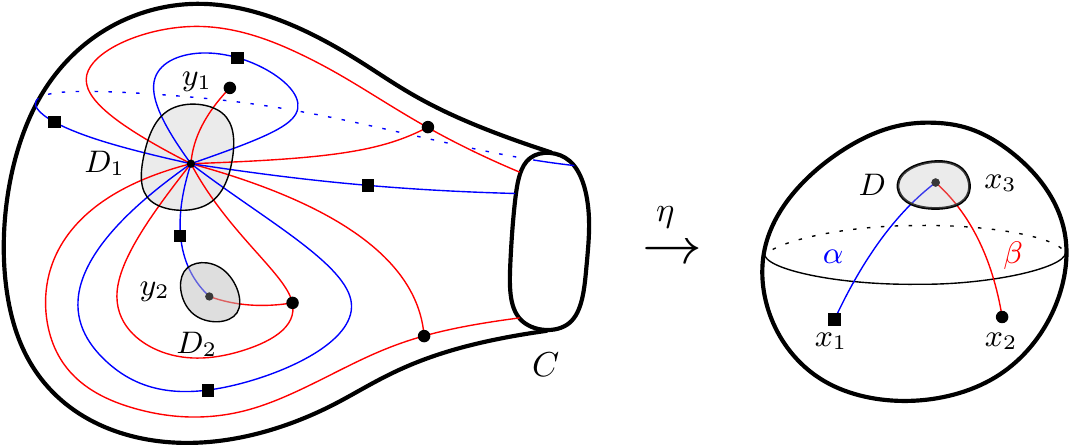}
\end{center}
\caption{ Case \textbf{(7.a)}. $\eta:D^2\rightarrow S^2$ defines a special almost orbifold  covering of degree $10$.}\label{fig:case7a}
\end{figure}

Let $\eta:D^2\rightarrow S^2$ be the map described in Fig.~\ref{fig:case7a} and let $\mathcal{O}'$ be the orbifold  with underlying surface $D^2$ and with two cone pints $y_1$  and $y_2$ of order $p_2=3$ and $p_3=7$ respectively, that is, $\mathcal{O}'=D^2(3,  7)$.   Then $\eta$  defines an almost orbifold covering of degree $10$ from $\mathcal{O}'$ to $\mathcal{O}$. 

The description of $\eta$ implies that  $\eta^{-1}(D)=D_1\sqcup D_2 \sqcup  C$ and that the restriction of $\eta$ to the disk $D_1$    defines an orbifold covering of degree $7$ from the orbifold $\mathcal{Q}'=D_1$ (a disk without cone points) onto  $\mathcal{Q}=D(p_3)$ (a disk with a single cone point $x_3$  of order $p_3=7$). Moreover,  $\eta|_{ D_2}:D_2\rightarrow D$  defines a  genuine   covering (of spaces) of degree one.   The exceptional boundary component $C$ is mapped  by $\eta$ onto $\partial D$ with degree $3$ which is smaller than $p_e=7$. Thus   $\eta$ is a special almost orbifold covering. 

It is not hard to see that    the restriction $\eta'$  of $\eta$ to the  surface  $$P:=D^2-\text{Int}(D_1)\cup \text{Int}(D_2)\subseteq D^2$$  
defines an orbifold covering of degree $10$ from the orbifold $\mathcal{O}_0'=P(p_2)$ (a pair of pants with a single cone point $y_1$ of order $p_2=3$) onto the orbifold $\mathcal{O}_0\subseteq \mathcal{O}$.   Moreover,  we have  $\eta_{\ast}'(\pi_1^o(\mathcal{O}_0'))= \langle  T_1, T_2, T_1T_2\rangle $ 
where $T_1= S_2$, 
$$T_2= S_1S_2S_1S_2S_1S_2^2 \cdot S_1S_2 \cdot S_2S_1 S_2^2S_1S_2^2S_1$$ 
and $T_3=(S_1S_2)^7$. By the injectivity of $\eta_{\ast}'$   we can assume that  $\pi_1^o(\mathcal{O}_0')=\langle T_1 , T_2 ,T_3\rangle.$
Using the relations in $G$  we see that  $i_{\ast} \circ \eta_{\ast}' (T_1)=s_2$, that 
\begin{eqnarray}
i_{\ast} \circ \eta_{\ast}' (T_2) & = &  s_1s_2s_1s_2s_1s_2^2 \cdot s_1s_2 \cdot s_2s_1 s_2^2s_1s_2^2s_1 \nonumber \\ 
          & = & s_3^{-3} s_1 \cdot  s_3^{-1} \cdot s_1  s_3^3 \nonumber
\end{eqnarray} 
and that  $i_{\ast} \circ \eta_{\ast}'(T_3)=(s_1s_2)^7=s_3^{-7}=1.$ The inclusion map $j:P\hookrightarrow D^2$ induces a surjective homomorphism $j_{\ast}:\pi_1^o(\mathcal{O}_0')\rightarrow \pi_1^o(\mathcal{O}')$   whose kernel is normally generated by $T_3$.   Consequently $\pi_1^o(\mathcal{O}')$ is generated by  
$$t_1:=j_{\ast}(T_1) \ \ \text{ and } \ \ t_2:=j_{\ast}(T_2).$$   
 As we clearly have $\eta_{\ast}\circ j_{\ast}=i_{\ast}\circ \eta_{\ast}'$ we see that $\eta_{\ast}$ sends the generating pair $(t_1, t_2^{\nu})$  of $\pi_1^o(\mathcal{O}')$   onto the pair $T'$.

\smallskip

For the reaming of this paper  $F$  will denote  the surface  obtained from a torus by removing an open  open disk   and   $\mathcal{O}'$  will denote  the orbifold with underlying surface $F$ and without cone points, that is, $\mathcal{O}'$ is isomorphic to the surface $F$.

\smallskip

\begin{figure}[h!]
\begin{center}
\includegraphics[scale=0.8]{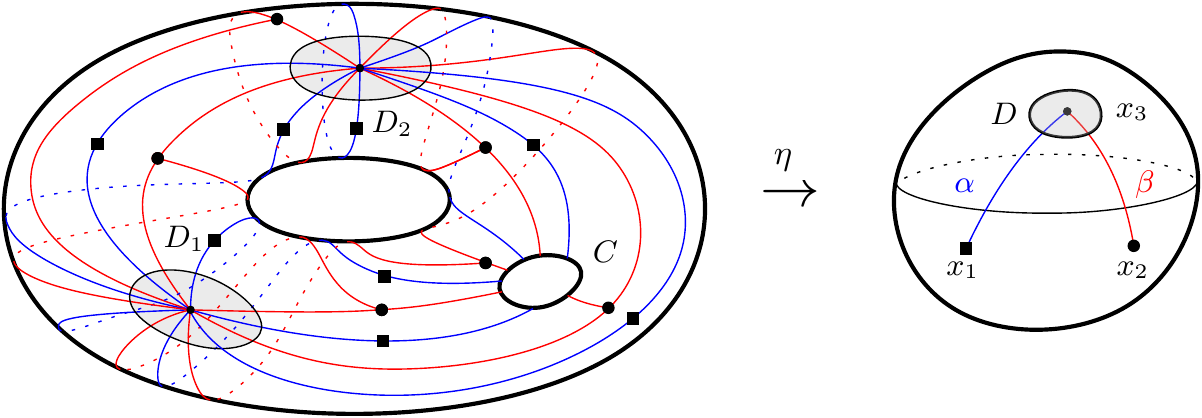}
\end{center}
\caption{ Case \textbf{(7.b)}. $\eta:F\rightarrow S^2$ defines a special  almost orbifold covering of degree $17$.}\label{fig:case7b}
\end{figure}

We now turn our attention to case (7.b). Thus we must find a special almost orbifold covering that represents  the pair 
$$T'= ((s_1s_2^2s_1s_2)^2s_2s_1s_2, (s_2(s_2s_1)^2)^2 s_2s_1).$$ 
The map $\eta:F\rightarrow S^2$  described  in Fig.~\ref{fig:case7b}  defines an almost orbifold covering of degree $17$ from the orbifold  $\mathcal{O}' $ to $\mathcal{O}$. The description of $\eta$ implies that
$$\eta^{-1}(D)=D_1\sqcup D_2 \sqcup  C$$  and that  
$$\eta|_{D_i}:D_i\rightarrow D$$  defines an orbifold covering of degree $7$ from the orbifold $\mathcal{Q}_i'=D_i$ (a disk without cone points) onto $\mathcal{Q}=D(p_3)$ (a disk with a single point of order $p_3=7$) for $i=1,2$. As   $\eta|_{C}:C\rightarrow \partial D$ has degree $4$ we conclude that   $\eta$ is a special almost  orbifold covering.
 
Put $F':=F-\text{Int}(D_1)\cup \text{Int}(D_2)$.  It is not hard to see that the  map   
$$\eta':=\eta|_{F'}: F'  \rightarrow D_0$$  
 defines an orbifold covering of degree $17$ from the  orbifold $\mathcal{O}_0'=F'$  onto the orbifold $\mathcal{O}_0\subseteq \mathcal{O}$.  Moreover, we have 
 $$\eta_{\ast}'(\pi_1^o(\mathcal{O}_0'))= \langle  T_1, T_2, T_3 , T_4, T_5\rangle \leq \pi_1^o(\mathcal{O}_0)$$ 
 where 
\begin{eqnarray}
T_1 & = & (S_1S_2^2S_1S_2)^2S_2S_1S_2  \nonumber \\
T_2 & = &  (S_2(S_2S_1)^2)^2 S_2S_1) \nonumber \\
T_3 & = &  S_1S_2^2\cdot (S_1S_2)^7\cdot S_1S_2^{2} \nonumber \\
T_4 & = & (S_1S_2)^{3}S_2S_1\cdot (S_1S_2)^7 S_1S_2 (S_1S_2)^3 \nonumber \\
T_5 & = &  (S_1S_2)^{7} \nonumber 
\end{eqnarray}
By the injectivity of $\eta_{\ast}'$   we can assume that   $\pi_1^o(\mathcal{O}_0')=\langle T_1 , T_2 ,T_3, T_4, T_5\rangle.$
Using the relations in $G$  we see that  $i_{\ast} \circ \eta_{\ast}' (T_1)= (s_1s_2^2s_1s_2)^2s_2s_1s_2$, that  $i_{\ast} \circ \eta_{\ast}' (T_2)= (s_2(s_2s_1)^2)^2 s_2s_1)$ and that $$i_{\ast} \circ \eta_{\ast}' (T_3)=i_{\ast} \circ \eta_{\ast}' (T_4)=i_{\ast} \circ \eta_{\ast}' (T_5)=1.$$

 Now the inclusion map    $j: F'\hookrightarrow F$  induces a surjective homomorphism  $$j_{\ast}:\pi_1^o(\mathcal{O}_0')\rightarrow \pi_1^o(\mathcal{O}')$$    whose kernel is normally generated by $\{T_3, T_4, T_5\}$.   Consequently $\pi_1^o(\mathcal{O}')$ is generated by  $$t_1:=j_{\ast}(T_1) \ \ \text{ and } \ \  t_2:=j_{\ast}(T_2).$$   
Since $\eta_{\ast}\circ j_{\ast}=i_{\ast}\circ \eta_{\ast}'$ we see that $\eta_{\ast}$ sends the  pair $(t_1, t_2)$   onto the pair $T$, and hence  the special almost orbifold covering $\eta:\mathcal{O}'\rightarrow \mathcal{O}$ represents  $T'$. 
 
\smallskip

\noindent{\textbf{Case (8).}}  Recall that  $p_1=p_2=3$ and $p_3\geqslant 4$ with $(p_3, 3)=1$. We must find a special almost orbifold covering that represents the generating pair $T'= (s_1s_2^2, s_2^2s_1).$  

\begin{figure}[h!]
\begin{center}
\includegraphics[scale=0.8]{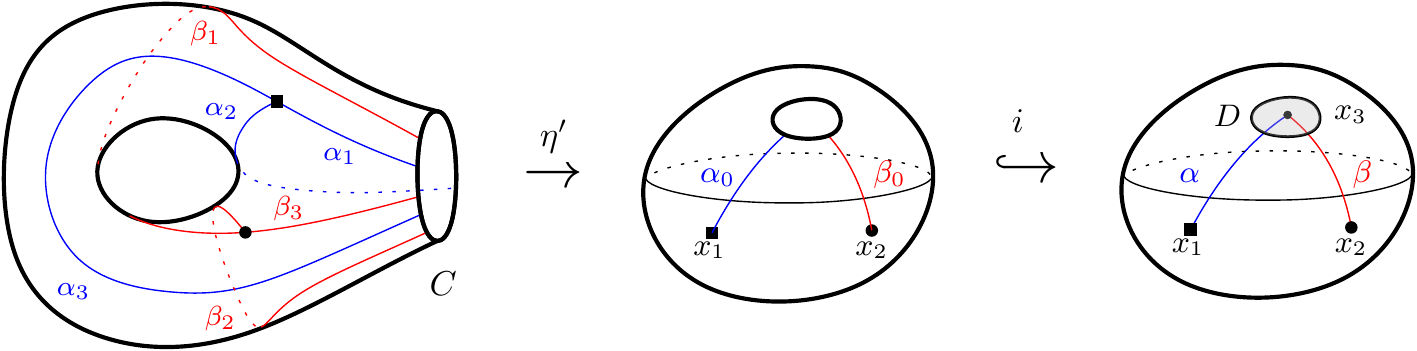}
\end{center}
\caption{ Case \textbf{(8)}.   $\eta =i\circ \eta' :F\rightarrow S^2$ defines an almost orbifold covering of degree $3$.}\label{fig:case9}
\end{figure}

Let $\eta':F\twoheadrightarrow D_0\subseteq S^2$ be the surjective map described in Fig.~\ref{fig:case9}, that is, $\eta'(\alpha_i)=\alpha_0$, $\eta'(\beta_i)=\beta$ for $1\leqslant i\leqslant 3$ and $\eta'$ maps each component of $$F- \alpha_1\cup\alpha_2\cup\alpha_3\cup \beta_1\cup \beta_2\cup\beta_3$$ 
homeomorphically onto $D_0-\alpha_0\cup\beta_0$. 

It is not hard to see that  $\eta'$  defines an  orbifold covering   of degree $3$ from $\mathcal{O}'$ onto  $\mathcal{O}_0$. Note further  that  the orbifold covering $\eta'$ corresponds  to the subgroup 
 $\langle  t_1, t_2\rangle$ of $\pi_1^o(\mathcal{O}_0)$ where $t_1= S_1S_2^2$ and $t_2= S_2^2S_1$. Thus   we  can  assume  that $\pi_1^o(\mathcal{O}')= \langle t_1 , t_2\rangle$.  Moreover,  $\eta'$  maps the boundary $C$ of $\mathcal{O}'$ onto  $\partial  D$  with degree $3$. 

It follows from Example~\ref{ex:almost} that $\eta:=i\circ \eta':F\rightarrow S^2$  defines an almost orbifold covering from $\mathcal{O}'$ to $\mathcal{O}$. Since $\eta|_C$ has degree $3$, which is smaller than $p_3$, we conclude that $\eta$ is a special almost orbifold covering. Now   the generating pair $(t_1, t_2)$ is clearly  mapped by $\eta_{\ast}=i_{\ast}\circ \eta_{\ast}'$ onto the pair $T'$.

\smallskip

\noindent{\textbf{Case (9).}} In this case   $p_1=2$, $p_2=4$ and $p_3\geqslant 7$ odd, and  $T'=( s_1s_2^2,  s_2^3  s_1  s_2^3)$.

\begin{figure}[h!]
\begin{center}
\includegraphics[scale=0.8]{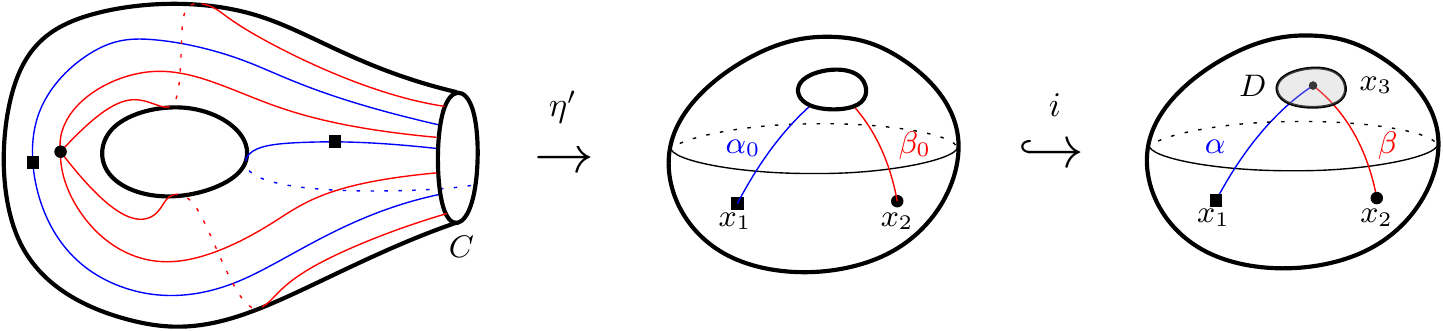}
\end{center}
\caption{ Case \textbf{(9)}.  $\eta =i\circ \eta ' :F\rightarrow S^2$ defines a special almost orbifold covering of degree $4$.}\label{fig:case10}
\end{figure}

Let $\eta':F\twoheadrightarrow D_0\subseteq S^2$ be the surjective map described in Fig.~\ref{fig:case10}. Thus  $\eta'$ defines an  orbifold covering of degree $4$ from $\mathcal{O}'$ onto  $\mathcal{O}_0$.  Moreover,  
$$(\eta_n')_{\ast}\pi_1^o(\mathcal{O}')=\langle T_1, T_2\rangle$$ 
where   $t_1= S_1S_2^2$  and $t_2= S_2^3  S_1  S_2^3$, consequently we may assume that $\pi_1^o(\mathcal{O}')$ is generated by the pair  $(t_1, t_2)$.  

It follows from Example~\ref{ex:almost} that $\eta:=i\circ \eta':F\rightarrow S^2$  defines an almost orbifold covering from $\mathcal{O}'$ to $\mathcal{O}$ which is special as  the degree of the map  $\eta|_C$  is equal to $4$ which  is smaller than $p_3$. The generating pair $(t_1, t_2)$  of $\pi_1^o(\mathcal{O}'$)   clearly gets mapped by $\eta_{\ast}=i_{\ast}\circ \eta_{\ast}'$ onto the pair $\tau'$.

\smallskip 

\noindent{\textbf{Case (10).}} In this case we have $p_1=2$, $p_2=3$ and $p_3\geqslant 7$ with $(p_3,6)=1$, and  $T'$ is equal to  $(s_1s_2s_1s_2^2, s_2^2s_1s_2s_1).$

\begin{figure}[h!]
\begin{center}
\includegraphics[scale=0.8]{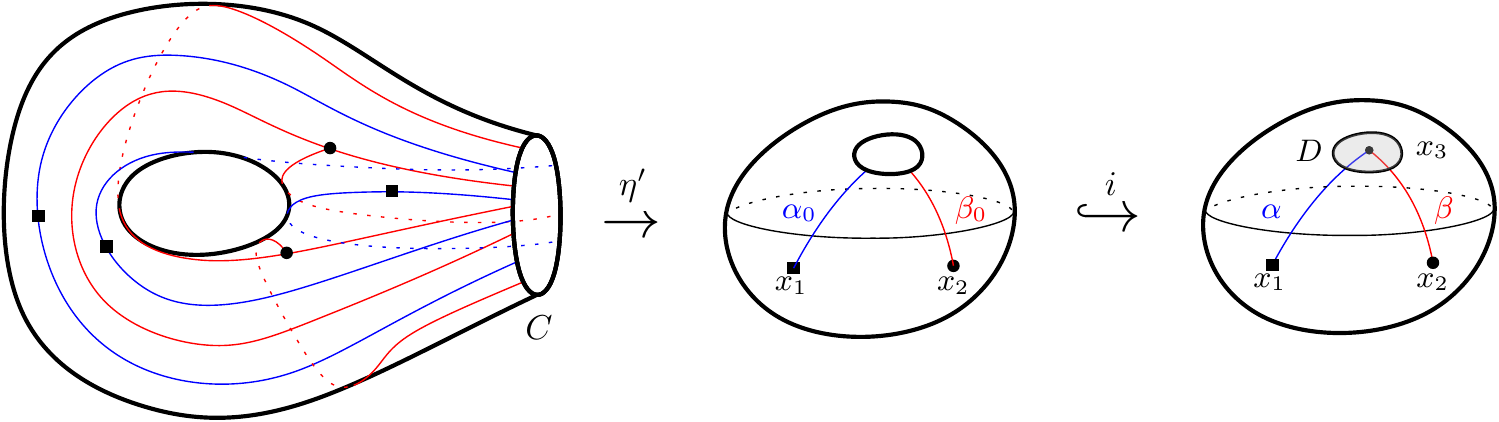}
\end{center}
\caption{ Case \textbf{(10)}.   $\eta =i\circ \eta ' :F\rightarrow S^2$ defines a special almost orbifold covering of degree $6$.}\label{fig:case8}
\end{figure}

The map $\eta':F\twoheadrightarrow D_0$  described   in Fig.~\ref{fig:case8} defines an orbifold covering of degree $6$ from $\mathcal{O}'$ onto $\mathcal{O}_0$ such that 
$$\eta_{\ast}'(\pi_1^o(\mathcal{O}'))  = \langle S_1S_2S_1S_2^2 ,  S_2^2S_1S_2S_1  \rangle.  $$
By applying  Example~\ref{ex:almost} we conclude that $\eta:=i\circ \eta'$ defines a almost orbifold covering from $\mathcal{O}'$ to $\mathcal{O}$. Moreover it is easy to see that 
$$ \eta_{\ast}(S_1S_2S_1S_2^2 ,  S_2^2S_1S_2S_1 )=T'$$  
which concludes the proof. 


\end{document}